\documentclass[]{amsart}
\usepackage[foot]{amsaddr}
\usepackage[a4paper,top=2.5cm,bottom=3cm,inner=3cm,outer=3cm]{geometry}

\usepackage[utf8]{inputenc}

\usepackage{amsmath,amsthm}
\usepackage{amsfonts,amssymb,mathrsfs}

\usepackage{xy}
\usepackage[bookmarks=false,breaklinks=true]{hyperref}
\usepackage{float}
\usepackage{graphicx}
\usepackage[percent]{overpic}
\usepackage{color}

\usepackage{epigraph}
\definecolor{myurlcolor}{rgb}{0,0,0.7}
\usepackage{hyperref}
\hypersetup{colorlinks,
linkcolor=myurlcolor,
citecolor=myurlcolor,
urlcolor=myurlcolor}
\usepackage[mathscr]{euscript}

\newcommand{\X}{{\mathbf X}}   

\newcommand{\C}{{\mathbb C}}  

\newtheorem*{thm*}{Theorem}

\begin{document}
\title{The Beauty of Roots}
\author[Baez]{John C.\ Baez} 
\address{Department of Mathematics, University of California, Riverside CA, 92521, USA}
\email{baez@math.ucr.edu}
\author[Christensen]{J.\ Daniel Christensen}
\address{University of Western Ontario, London, Ontario, Canada}
\email{jdc@uwo.ca}
\author[Derbyshire]{Sam Derbyshire}
\address{16 Place Saint Clément, Luz St Sauveur, France}
\email{sam@well-typed.com}
\date{October 1, 2023}
\maketitle

\begin{center}
\includegraphics[width=6in]{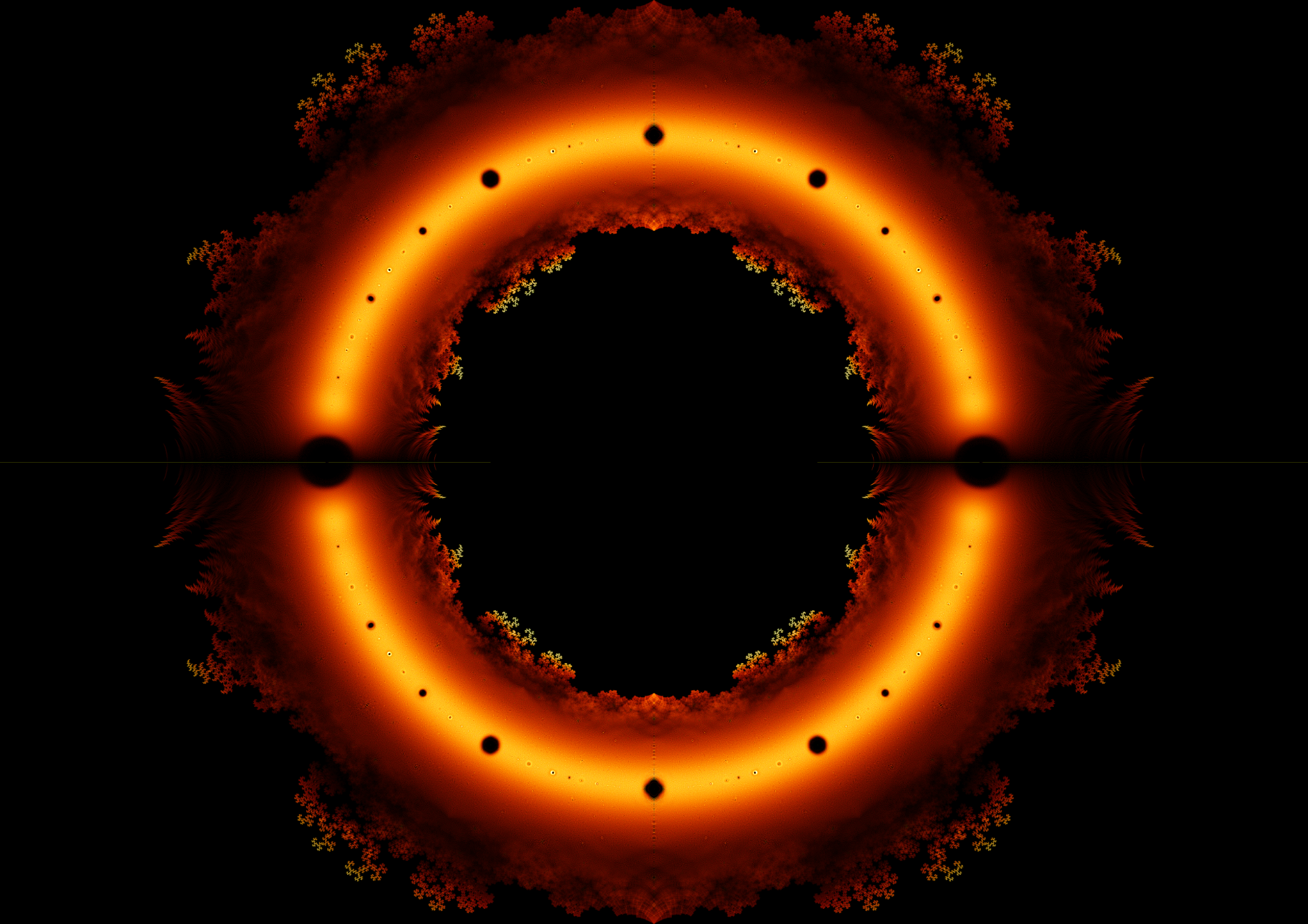}
\vskip 0.5em
Figure 1. Roots of all polynomials of degree 23 whose coefficients are $\pm 1$.
The brightness shows the number of roots per pixel.  
\end{center}
\vskip 1em

One of the charms of mathematics is that simple rules can generate
complex and fascinating patterns, which raise questions whose answers require
profound thought.   For example, if we plot the roots of all polynomials of degree $23$
whose coefficients are all $1$ or $-1$, we get an astounding picture, shown in Figure 1.

More generally, define a \textbf{Littlewood polynomial} to be a polynomial $p(z) = \sum_{i=0}^d a_i z^i$ with each coefficient $a_i$ equal to $1$ or $-1$. Let $\X_n$ be the set of complex numbers that are roots of some Littlewood polynomial with $n$ nonzero terms (and thus degree $n-1$).   The 4-fold symmetry of Figure 1 comes from the fact that if $z \in \X_n$ so are $-z$ and $\overline{z}$.  The set $\X_n$ is also invariant under the map $z \mapsto 1/z$, since if $z$ is the root of some Littlewood polynomial then $1/z$ is a root of the polynomial with coefficients listed in the reverse order.  

It turns out to be easier to study the set
\[        \X = \bigcup_{n = 1}^\infty \X_n = 
\{ z \in \C | \; z \; \textrm{is the root of some Littlewood 
polynomial} \} .\]
If $n$ divides $m$ then $\X_n \subseteq \X_{m}$, so $\X_n$ for a highly divisible
number $n$ can serve as an approximation to $\X$, and this is why we drew $\X_{24}$.

Some general properties of $\X$ are understood.   It is easy to show that $\X$ is 
contained in the annulus $1/2 < |z| < 2$.  On the other hand, Thierry Bousch 
showed \cite{Bou1} that the closure of $\X$ 
contains the annulus $2^{-1/4} \le |z| \le 2^{1/4}$.   This means that the holes 
near roots of unity visible in the sets $\X_d$ must eventually fill in as we take the 
union over all degrees $d$.    More surprisingly, Bousch showed in 1993 that 
the closure $\overline{\X}$ is connected and locally path-connected \cite{Bou3}.  
It is worth comparing the work of Odlyzko and Poonen \cite{OP}, who previously
showed similar result for roots of polynomials  whose coefficients are all $0$ or $1$.

\begin{center}
\includegraphics[width=3in]{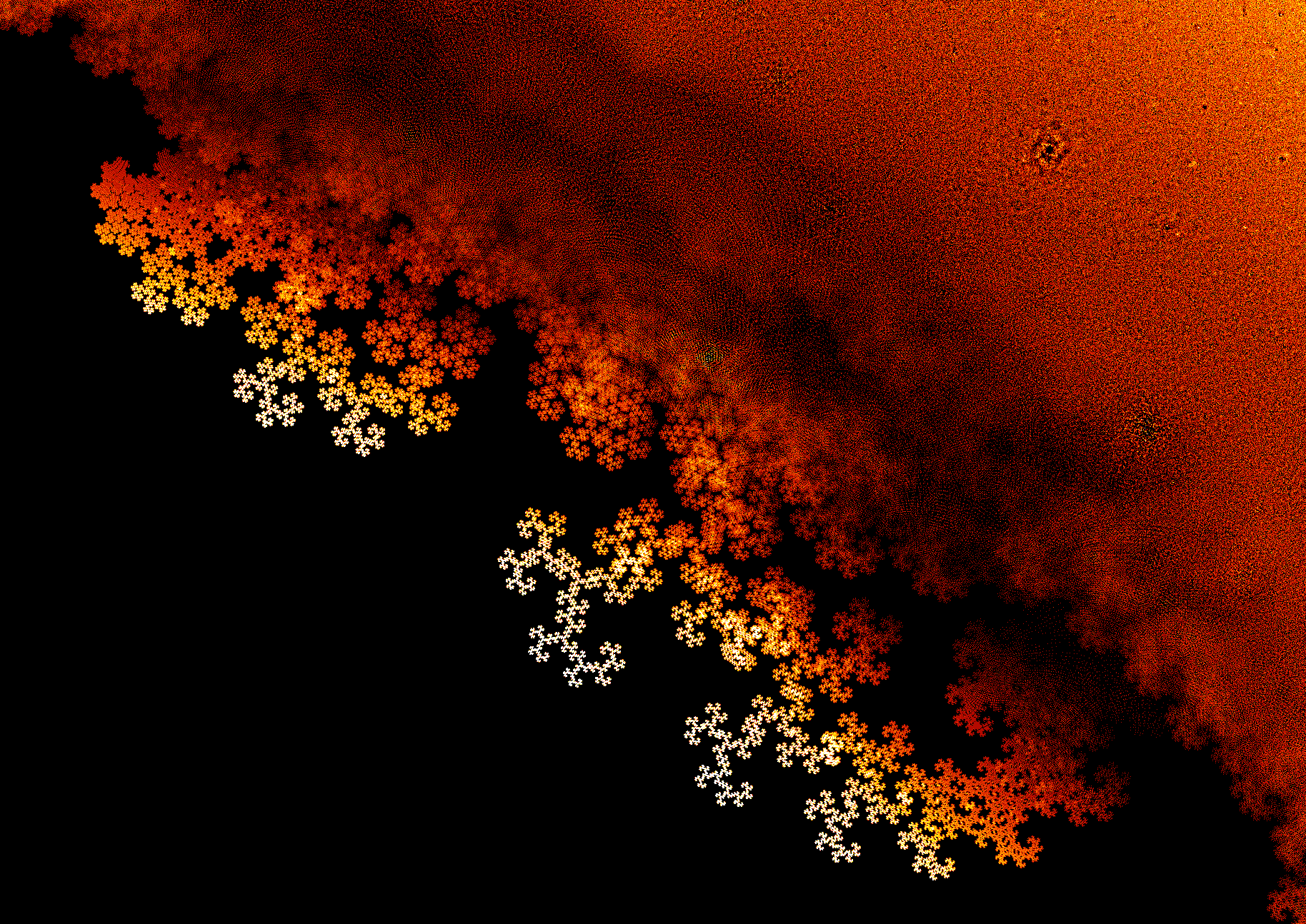}
\vskip 0.5em
Figure 2. The region of $\X_{24}$ near the point $z = \frac{1}{2}e^{i/5}$.
\end{center}

The big challenge is to understand the diverse, complicated and beautiful patterns that appear in different regions of the set $\X$.   There are websites that let you explore 
and zoom into this set online \cite{C,Egan,V}.  Different regions raise different questions.

For example, what is creating the fractal patterns in Figure 2 and elsewhere?   An anonymous contributor suggested a fascinating line of attack which was further developed by Greg Egan \cite{Egan}. Define two functions from the complex plane to itself, depending on a complex parameter $q$: 
\[  
 f_{+q}(z) = 1 + q z , \qquad  f_{-q}(z) = 1 - q z .
\]
When $|q| < 1$ these are both contraction mappings, so by a theorem
of Hutchinson \cite{H} there is a unique nonempty compact 
set $D_q \subseteq \C$ with
\[         D_q = f_{+q}(D_q) \cup f_{-q}(D_q)  .\]
We call this set a \textbf{dragon}, or the {\boldmath{$q$\textbf{-dragon}}} to
be specific.   And it seems that \textit{for $|q| < 1$, the portion of the set $\X$ in a small neighborhood of the point $q$ tends to look like a rotated version of $D_q$}.  

Figure 3 shows some examples.  To precisely describe what is going on, much less prove it, would take real work.   We invite the reader to try.  A heuristic explanation is known, which can serve as a starting point \cite{Baez,Egan}.   Bousch \cite{Bou3} has also proved this related result:

\begin{thm*} \label{dragon.1}
For $q \in \C$ with $|q| < 1$, we have $q \in
\overline{\X}$ if and only if $0 \in D_q$.  When this holds, the set $D_q$ is connected.
\end{thm*}

\begin{center}
\includegraphics[height=2in]{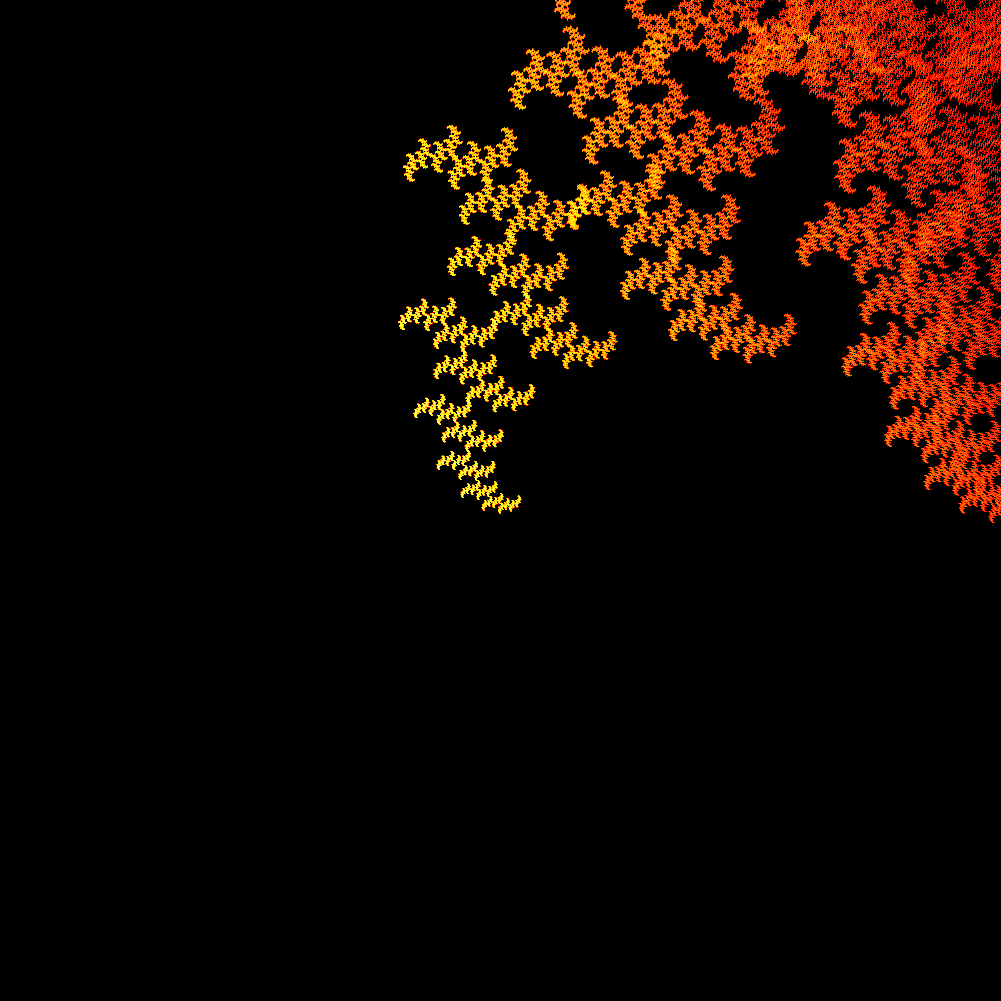} \quad
\includegraphics[height=2in]{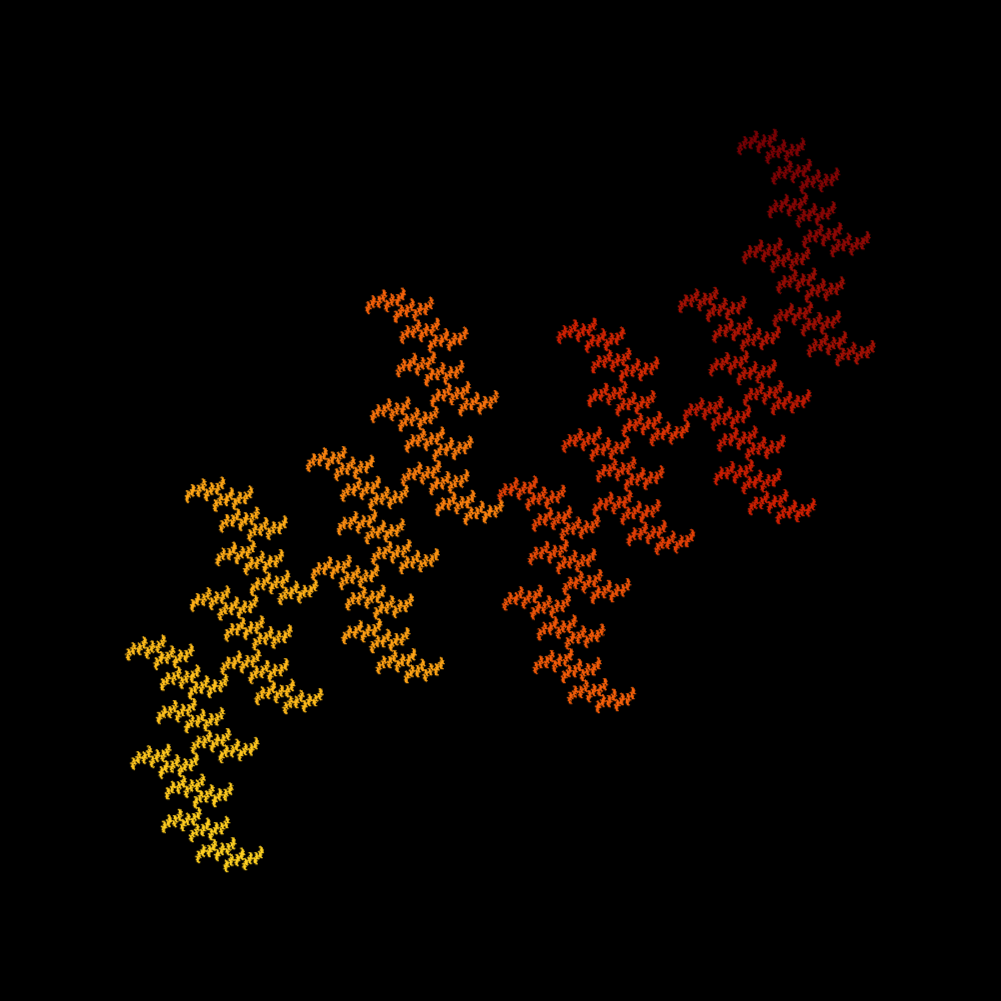} \break \vskip 0.1em
\includegraphics[height=2in]{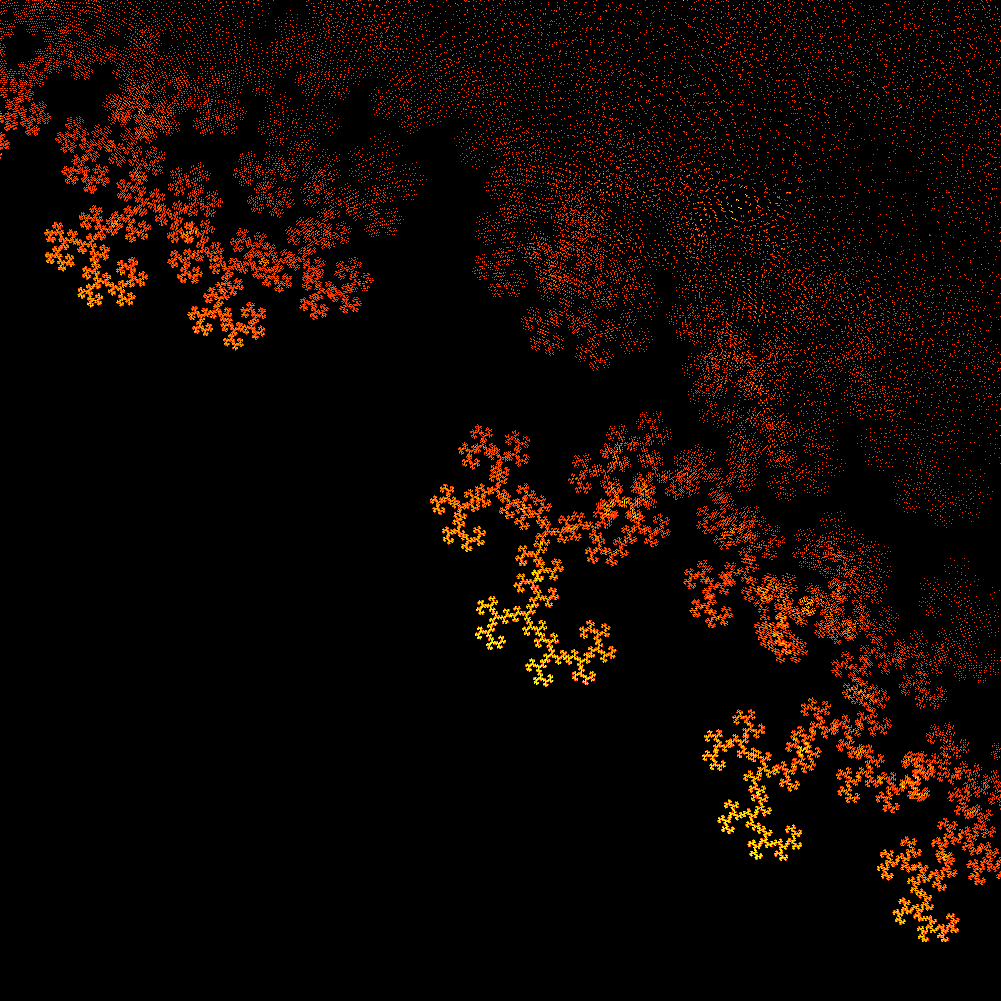} \quad
\includegraphics[height=2in]{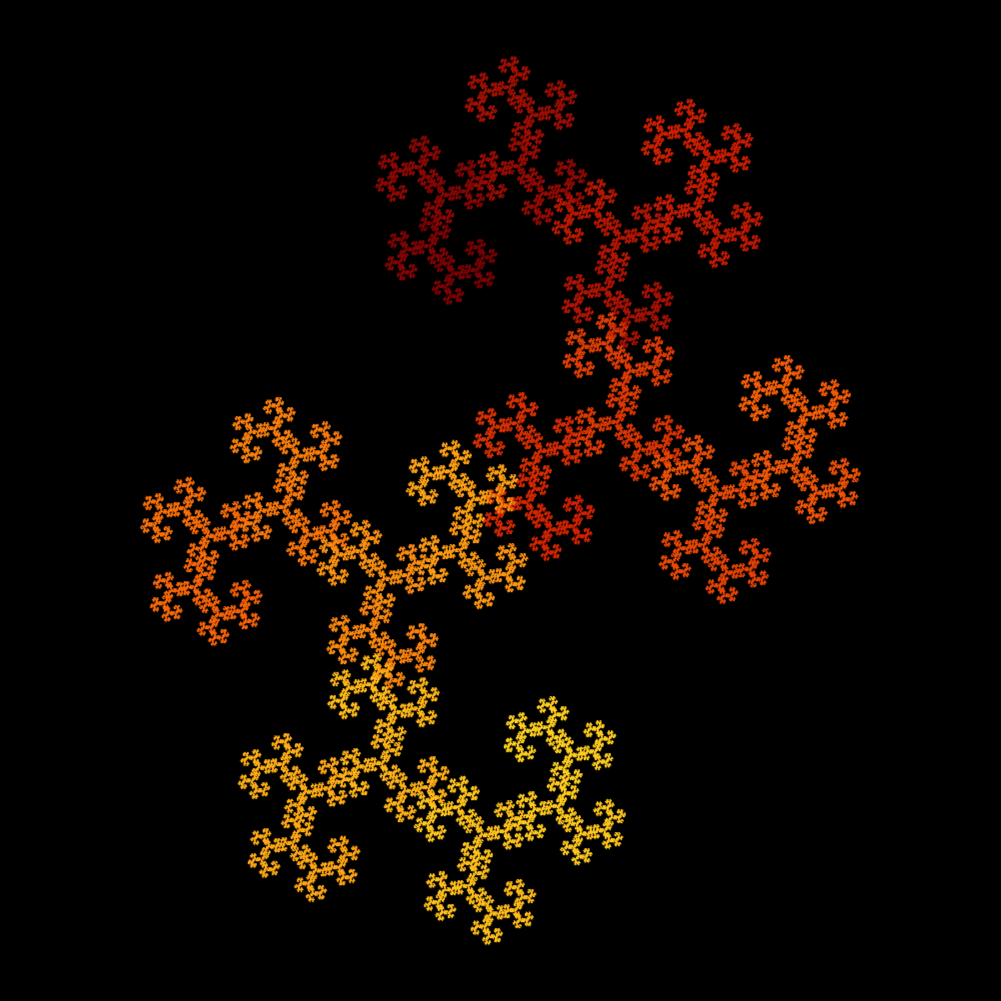}
\vskip 0.5 em
Figure 3.  
Top: the set $\X$ near $q = 0.594 + 0.254i$ at left, and the set $D_q$ at right. \\
Bottom: the set $\X$ near $q = 0.375453 + 0.544825i$ at left, and the set $D_q$ at right.
\end{center}

\end{document}